\documentclass{article}

\usepackage{amssymb, amsmath}

  \newtheorem{lemma}{Lemma}%

\newtheorem{remark}{Remark}

\newcommand{\vep}{\varepsilon}

\newcommand{\QED}{\hfill
$\underline{\underline{QED}}$}
\newcommand{\ds}{\displaystyle}

\newcommand{\la}{\lambda}

\newcommand{\beq}{\begin{equation}}
\newcommand{\eeq}{\end{equation}}
\newcommand{\beqn}{\begin{eqnarray}}
\newcommand{\eeqn}{\end{eqnarray}}
\newcommand{\beqnn}{\begin{eqnarray*}}
\newcommand{\eeqnn}{\end{eqnarray*}}
\newcommand{\nn}{|\!|}

\newcommand{\rr}{\mathbb{R}}

\begin{document}
\title{Landesman-Laser Conditions and Quasilinear Elliptic Problems}

\author{Nikolaos Halidias \\
University of the Aegean \\ Department of Statistics and Actuarial
Science
\\ Karlovassi, 83200 \\ Samos \\ Greece \\ email: nick@aegean.gr}

\maketitle


\begin{abstract} In this paper we consider two elliptic problems. The first one
is a Dirichlet problem while the second is Neumann. We extend all
the  known results concerning  Landesman-Laser conditions by using
the Mountain-Pass theorem with the Cerami $(PS)$
condition.\footnote{2000 Mathematics Subject Calssification:
35J25, 35J60
 }
\end{abstract}

\section{Introduction}

In this paper using the well-known Mountain-Pass Theorem and
Cerami (PS) (see \cite{BBF}) condition we extend all the known
results concerning quasilinear elliptic problems at resonance
satisfying the Landesman-Laser conditions.

 Before we proceed we must state some
well-known definitions and facts. Let $X$ be a Banach space. We
say that a functional $I:X \to \rr$ satisfies the $(PS)_c$
condition if for any sequence such that $|I(u_n)|\leq M$ and
$(1+\nn u_n \nn)<I^{'}(u_n),\phi> \to 0$ for all $\phi \in X$ we
can show that there exists a convergent subsequence.

 Consider the first eigenvalue $\lambda_1$ of
$(-\Delta_p,W^{1,p}_o(\Omega))$. From Lindqvist \cite{Li}  we know
that $\lambda_1
> 0 $ is isolated and simple, that is any two solutions $u,v$ of

\begin{equation}
\left\{
\begin{array}{l}
-\Delta_pu= -div(\nn Du \nn^{p-2}Du)=\lambda_1|u|^{p-2}u  \mbox{
a.e. on }\Omega \ \\ u\mid_{\partial \Omega} = 0, 2 \leq p<\infty
\end{array}
\right\}
\end{equation}

satisfy $u=cv$ for some $c \in \rr$. In addition, the
$\lambda_1$-eigenfunctions do not change sign in $\Omega$. Finally
we have the following variational characterization of $\lambda_1$
(Rayleigh quotient): $$ \lambda_1 = \inf \bigg[\frac{\nn Du
\nn^p_p}{\nn u \nn^p_p}:u \in W^{1,p}_o(\Omega),u \neq 0 \bigg] $$

We can define
\begin{eqnarray*}
\hat{\la}_2 = \inf \{ \la > 0: \la \mbox{ is an eigenvalue of
}(-\Delta_p, W^{1,p}_o(\Omega)), \la \neq \la_1 \} > \la_1.
\end{eqnarray*}
Anane-Tsouli \cite{AT} had proved that the second eigenvalue,
$\la_2$, is equal with $\hat{\la}_2$ and has a variational
characterization.

Let us state our first problem. Let $\Omega \subseteq \rr^n$ be a
bounded domain with smooth enough boundary $\partial \Omega$. The
Dirichlet problem is
\begin{equation}
\left\{
\begin{array}{lll}
- div \bigl(\nn Du(x) \nn^{p-2} Du(x)\bigr) -
\la_1|u(x)|^{p-2}u(x) = f(x,u(x)) \mbox{ a.e. on } \Omega
\\ u = 0
 \mbox{ a.e. on } \partial \Omega, \;\; 2 \leq p < \infty.
\end{array}
\right.
\end{equation}

Recently Bouchala-Drabek \cite{BD} had considered the above
problem and they derive via the Saddle-Point Theorem a weak
solution when the right-hand side satisfies an extended type of
Landesman - Laser conditions. Here we extended more that
conditions by using the following hypotheses.

 $ H(f):$ $f:\Omega \times    \rr \to \rr$ is a
Carath\'eodory function. Moreover,
\begin{enumerate}
\item [(i)] for almost all $x \in \Omega$ and all $u \in \rr$, $|f(x,u)|
\leq a(x) + c_1|u|^{p-1} $, $a(x) \in L^{\infty}(\Omega)$;
\item [(ii)] uniformly for all $x \in \Omega$ we have
 $\ds{\limsup_{u\to 0}} \frac{p F(x,u)}{|u|^p} \leq
\theta (x) \leq 0$  with $\theta (x) \in L^{\infty}(\Omega)$ and
$\int_{\Omega} \theta (x)) |u_1(x)|^pdx < 0$, $u_1$ is the first
eigenvalue and  $F(x,u) = \int_o^{u}f(x,r)dr$;
\item[(iii)]uniformly for almost all $ x \in \Omega$ we have that
\begin{eqnarray*}
\lim_{|u| \to \infty} \frac{F(x,u)}{|u|^p} = 0,
\end{eqnarray*}
moreover, there exists  a function $h : \rr^{+} \to \rr^{+}$ with
the property $\liminf \frac{h(a_nb_n)}{h(b_n)} \geq 1$, $h(b_n)
\to \infty$ when $a_n \to a>0$ and $b_n \to +\infty$ and another
function $\mu (x) \in L^{\infty}(\Omega)$ with $\int_{\Omega} \mu
(x) dx >0$, such that
$$\liminf_{|u| \to \infty} \frac{p F(x,u) - f(x,u)u}{h(|u|)} \geq \mu (x).$$
\end{enumerate}

In order to use the Mountain-Pass Theorem we must define the
energy functional of our problem. Let $I : X \to \rr$ be such that
$I(u) = \frac{1}{p} \nn Du \nn_p^p - \frac{\la_1}{p} \nn u \nn_p^p
- \int_{\Omega} F(x,u(x))dx$. It is well-known that $I$ is a $C^1$
functional and its critical points are in fact weak solutions to
problem (2).

\begin{lemma}
If hypotheses $H(f)(i),(iii)$ holds, then the energy functional
satisfies the $(PS)_c$ condition.
\end{lemma}

{\bf Proof:} Let $X = W^{1,p}_o(\Omega)$.  Suppose that there
exists a sequence $\{ u_n \} \subseteq X$ such that $|I(u_n)| \leq
M$ and
\begin{eqnarray}
<I^{'}(u_n),\phi> \leq \vep_n \frac{\nn \phi \nn_{1,p}}{1 + \nn
u_n \nn_{1,p}}.
\end{eqnarray}
Suppose that $\nn u_n \nn_{1,p} \to \infty$. Let $y_n(x) =
\frac{u_n(x)}{\nn u_n \nn_{1,p}}$.

From the first inequality we have
\begin{eqnarray}
| \frac{1}{p} \nn Du_n \nn_p^p  - \frac{\la_1}{p} \nn u_n \nn_p^p
- \int_{\Omega} F(x,u_n(x))dx | \leq M.
\end{eqnarray}

From $H(f)(iii)$ we know that $\lim_{|u| \to \infty}
\frac{F(x,u)}{|u|^p} = 0$. It is easy to see that also $\lim_{n
\to \infty} \int_{\Omega}\frac{F(x,u(x))}{\nn u \nn^p_{1,p}} dx
=0$.

Dividing this inequality with $\nn u_n \nn_{1,p}^p$ and using
$H(f)(iii)$ we  arrive to the conclusion that $\nn Dy \nn_p \leq
\nn Dy_n \nn_p \to \la_1 \nn y \nn_p$. So using the uniform
convexity we arrive at the conclusion that $y_n \to y$ strongly in
$X$  and that $y = u_1 (x)$. Note that $\nn y_n \nn_{1,p} = 1$. So
we can say that $|u_n(x)| \to \infty$.

Choosing now $\phi = u_n$ in (3) and substituting with (4) we
arrive at
\begin{eqnarray*}
-M - \vep_n \frac{\nn u_n \nn_{1,p}}{1 + \nn u_n \nn_{1,p}} \leq
\int_{\Omega} (p F(x,u_n(x)) - f(x,u_n(x))u_n(x))dx \leq M +\vep_n
\frac{\nn u_n \nn_{1,p}}{1 + \nn u_n \nn_{1,p}}.
\end{eqnarray*}
Dividing now the last inequality with $h(\nn u_n \nn_{1,p})$ we
obtain
\begin{eqnarray*}
\frac{-M - \vep_n \frac{\nn u_n \nn_{1,p}}{1 + \nn u_n
\nn_{1,p}}}{h(\nn u_n \nn_{1,p})} \leq \int_{\Omega} \frac{p
F(x,u_n(x)) - f(x,u_n(x))u_n(x)}{h(|u_n(x)|)} \frac{h(|y_n(x)| \nn
u_n \nn_{1,p})}{h(\nn u_n \nn_{1,p})} dx \leq  \\ \frac{M + \vep_n
\frac{\nn u_n \nn_{1,p}}{1 + \nn u_n \nn_{1,p}}}{h(\nn u_n
\nn_{1,p})}.
\end{eqnarray*}
From this we can see that
\begin{eqnarray*}
\liminf_{n \to \infty}\int_{\Omega} \frac{p F(x,u_n(x)) -
f(x,u_n(x))u_n(x)}{h(|u_n(x)|)} \frac{h(|y_n(x)| \nn u_n
\nn_{1,p})}{h(\nn u_n \nn_{1,p})} dx \leq 0.
\end{eqnarray*}
Using lemma Fatou and  $H(f)(iii)$ we obtain the contradiction.
That is $u_n$ is bounded. Using well-known arguments we can also
show that in fact $u_n$ has a strongly convergent subsequence (see
\cite{BD}).

\QED

\begin{lemma}
If $H(f)$ holds, then there exists some $\rho,a > 0$ such that for
all $u \in W^{1,p}_o(\Omega)$ with $\nn u \nn_{1,p} = \rho$ we
have $I(u) \geq a
>0$.
\end{lemma}

{\bf Proof:} We shall show  that there exists $\rho > 0$ such that
$I(u) \geq a > 0$ with $\nn u \nn_{1,p} = \rho$. To this end,
 we'll show that
for every sequence $\{ u_n \}_{n \geq 1} \subseteq W^{1,p}_o$ with
$\nn u_n \nn = \rho_n \to 0$ we have $I(u_n) \downarrow 0$.
Suppose that it is not true. Then there exists a sequence as above
such that $I(u_n) \leq 0.$ Since $\nn u_n \nn_{1,p} \to 0$ we have
 $u_n(x) \to 0$ a.e. on $\Omega$.

So we have
\begin{eqnarray}
\nn Du_n \nn_p^p -\la_1 \nn u_n \nn_p^p \leq \int_{\Omega} p
F(x,u_n(x))dx.
\end{eqnarray}

Let $y_n (x) = \frac{u_n(x)}{\nn u_n \nn_{1,p}}$. Also, from
$H(f)(ii)$ we have uniformly for all $x \in \Omega$ that for all
$\vep
> 0$ we can find
 $\delta > 0 $ such that for $|u| \leq \delta$ we have
\begin{eqnarray*}
p F(x,u) \leq \theta (x) |u|^p + \vep |u|^p.
\end{eqnarray*}

On the other hand from hypothesis $H(f)(i)$  we have that there
exists some $c_1,c_2$ such that $p F (x,u) \leq c_1|u|^p +
c_2|u|^{p^*}+p|u|$ for almost all $x \in \Omega$ and all $u \in
\rr$. Thus we can always find some $\gamma > 0$ such that $p
F(x,u) \leq (\theta(x)+\vep)|u|^p + \gamma |u|^{p^*}$. Indeed,
choose $\gamma \geq |c_1-\theta (x)-\vep||\delta|^{p-p^*} +
c_2+p|\delta|^{1-p^*}$.

Then we obtain,
\begin{eqnarray}
\nn Du_n \nn_p^p -\la_1 \nn u_n \nn_p^p \leq  \int_{\Omega}
(\theta (x) + \vep) |u_n(x)|^pdx + \gamma \int_{\Omega}
|u_n(x)|^{p^*}dx.
\end{eqnarray}

 Dividing  inequality  (6) with $\nn u_n
\nn^p_{1,p}$, we have
\begin{eqnarray*}
 \nn Dy_n \nn^p  -\la_1 \nn y_n
\nn^p_p \leq
 \int_{\Omega} (\theta (x) + \vep) |y_n(x)|^pdx +
 \gamma \frac{\int_{\Omega} |u_n(x)|^{p^*}dx}{\nn u_n \nn^p_{1,p}}\leq \\ \vep \nn
y_n \nn^p_p +\gamma_1 \nn u_n \nn^{p^*-p}_{1,p},
\end{eqnarray*}
 recall that $W^{1,p}_o(\Omega)$ is continuously embedded on
 $L^{p^*}(\Omega)$.

 Using the variational characterization of the first
eigenvalue we have that
\begin{eqnarray}
 0
  \leq \nn Dy_n \nn^p_p -\la_1 \nn y_n \nn^p_p \leq \int_{\Omega} \vep
  |y_n(x)|^pdx+\gamma_1 \nn u_n \nn^{p^*-p}_{1,p} .
\end{eqnarray}

Recall that $\nn y_n \nn = 1$ so $y_n \to y$ weakly in
$W^{1,p}_o(\Omega)$, $y_n(x) \to y(x)$ a.e. on $\Omega$. Thus,
from  inequality (7) we have that $\nn Dy_n \nn \to \la_1 \nn y
\nn$. Also, from the weak lower semicontinuity of the norm we have
$\nn Dy \nn \leq \liminf \nn Dy_n \nn \to \la_1 \nn y \nn$. Using
the variational characterization of the first eigenvalue we have
that $\nn Dy \nn = \la_1 \nn y \nn$. Recall that $y_n \to y$
weakly in $W^{1,p}_o(\Omega)$ and $\nn Dy_n \nn \to \nn Dy \nn$.
So, from the Kadec-Klee property we obtain $y_n \to y$ in
$W^{1,p}_o(\Omega)$ and since $\nn y_n \nn =1$ we have that $\nn y
\nn =1$. That is, $y \neq 0$ and from the equality $\nn Dy \nn =
\la_1 \nn y \nn$ we have that $y(x) =u_1(x)$.

Dividing now (6) with $\nn u_n \nn^p_{1,p}$ and using the
variational characterization of the first eigenvalue we have, that
for every $\vep > 0$ there exists some $n_o$ such that for $n \geq
n_o$ we have
\begin{eqnarray*}
0 \leq  \int_{\Omega} (\theta (x) +\vep) |y_n(x)|^pdx + \gamma_1
\nn u_n \nn^{p^*-p}_{1,p}.
\end{eqnarray*}

So in the limit we obtain
\begin{eqnarray*}
\vep \nn u_1 \nn_p^p \geq \int_{\Omega} (-\theta (x) ) |u_1(x)|^p
dx \mbox{ for every }\vep >0.
\end{eqnarray*}

 So this is a contradiction. So there exists $\rho > 0$ such
that $I(u) \geq a >  0$ for all $u \in W^{1,p}_o(\Omega)$ with
$\nn u \nn_{1,p} = \rho$.

\QED

\begin{lemma}If hypotheses $H(f)$ holds, then there exists some $e
\in W^{1,p}_o(\Omega)$ with $I(e) \leq 0$.
\end{lemma}

{\bf Proof:} We will show that there exists some $a \in \rr$ such
that $I(a|u_1|) \leq 0$. Suppose that this is not the case. Then
there exists a sequence $a_n \in \rr$ with $a_n \to \infty$ and
$I(a_n |u_1|) \geq c
>0$.

We can easily see that
\begin{eqnarray*}
(-\frac{F(x,u)}{u^p})^{'} = \frac{p F(x,u)- f(x,u) u}{u^{p+1}} =
\frac{p F(x,u) - f(x,u) u}{h(|u|)} \frac{h(|u|)}{u^{p+1}} \\ \geq
(\mu
 (x)-\vep) \frac{1}{u^{p+1}} = \frac{\mu (x)-\vep}{p}
 (-\frac{1}{u^p})^{'},
 \end{eqnarray*}
 for big enough $u \in \rr$.

We can say then
\begin{eqnarray*}
\int_t^s (-\frac{F(x,u)}{u^p})^{'}du \geq \int_t^s \frac{\mu
(x)-\vep}{p}
 (-\frac{1}{u^p})^{'}du.
 \end{eqnarray*}
Take now $s \to \infty$ and using $H(f)(iii)$ we obtain
\begin{eqnarray*}
F(x,t) \geq \frac{\mu (x)}{p},
\end{eqnarray*}
for big enough $t \in \rr$.
 From this we obtain
\begin{eqnarray*}
\limsup_{a_n \to \infty} I(a_n |u_1|) \geq \liminf_{a_n \to
\infty} I(a_n |u_1|) \geq 0 \Rightarrow \\
\limsup_{a_n \to \infty} \int_{\Omega} -F(x,a_n |u_1(x)|)dx \geq 0
\Rightarrow \\
\int_{\Omega} \frac{-\mu (x)}{p}dx \geq 0.
\end{eqnarray*}

 Then using $H(f)(iii)$ we obtain the
contradiction.

\QED

The existence of the nontrivial solution follows from the
Mountain-Pass Theorem.

\section{Neumann Problems}
Let $X = W^{1,p}(\Omega)$. Before we start let us mention some
facts. It is well-known that $W^{1,p}(\Omega) = \rr\oplus W$ with
$W = \{ u \in X: \int_{\Omega} u(x)dx = 0 \}$. We introduce the
following number,
\begin{eqnarray*}
 \la_1 = \inf \{ \frac{\nn Dw \nn_p^p}{\nn w \nn_p^p}: w \in W, w
 \neq 0 \}.
 \end{eqnarray*}

 From Papalini \cite{Papa} we know that $\la_1 >0$ and if $w
 \in W$ is such that $\nn w \nn_p  = 1, \nn Dw \nn_p = \la_1$ then
 $w$ is an eigenunction of the following problem,
\begin{equation}
\left\{
\begin{array}{lll}
- div \bigl(\nn Du(x) \nn^{p-2} Du(x)\bigr) = \la_1
|u(x)|^{p-2}u(x) \mbox{ a.e. on } \Omega
\\ -\frac{\partial u}{\partial n_p} = 0
 \mbox{ a.e. on } \partial \Omega, \;\; 2 \leq p < \infty.
\end{array}
\right.
\end{equation}
Now we are ready to state our second problem.

The problem under consideration is the following:
\begin{equation}
\left\{
\begin{array}{lll}
- div \bigl(\nn Du(x) \nn^{p-2} Du(x)\bigr) = f(x,u(x)) \mbox{
a.e. on } \Omega
\\ -\frac{\partial u}{\partial n_p} = g(x,u(x))
 \mbox{ a.e. on } \partial \Omega, \;\; 2 \leq p < \infty.
\end{array}
\right.
\end{equation}

Let us state the hypotheses on the data. Set $F(x,u) = \int_o^{u}
f(x,r)dr, \quad G(x,u) = \int_o^{u} g(x,r)dr$.

$ H(f,g):$ $f,g:\Omega \times    \rr \to \rr$ are Carath\'eodory
functions. Moreover,
\begin{enumerate}
\item [(i)] for almost all $x \in \Omega$ and all $u \in \rr$, $|f(x,u)|
\leq a(x) + c_1|u|^{p-1} $, $|g(x,u)| \leq a(x) + c_1|u|^{p-1}$,
$a(x) \in L^{\infty}(\Omega)$;
\item [(ii)] uniformly for all $x \in \Omega$ we have
 $\ds{\limsup_{u\to 0}} \frac{p F(x,u)}{|u|^p} \leq
\theta (x) \leq \la_1$  with $\theta (x) \in L^{\infty}(\Omega)$
and $\int_{\Omega}(\la_1- \theta (x)) |w(x)|^pdx > 0$, $\lim_{|u|
\to 0} \frac{G(x,u)}{|u|^p} = 0$ for any $w \in W$ an
eigenfunction to $\la_1$.
\end{enumerate}

Finally, we have the following hypothesis,

H(fg): uniformly for almost all $ x \in \Omega$ we have that
\begin{eqnarray*}
\lim_{|u| \to \infty} \frac{F(x,u)}{|u|^p} = 0, \quad \lim_{|u|
\to \infty}\frac{G(x,u)}{|u|^p} = 0,
\end{eqnarray*}
and suppose that there is a function $h : \rr^{+} \to \rr^{+}$
with the property $\liminf \frac{h(a_nb_n)}{h(b_n)} \geq 1$,
$h(b_n) \to \infty$ when $a_n \to a>0$ and $b_n \to +\infty$ such
that
\begin{eqnarray*}
& & \liminf_{|u| \to \infty} \frac{p F(x,u)-f(x,u)u}{h(|u|)} \geq
\mu (x), \\
& & \liminf_{|u| \to \infty} -\frac{p G(x,u)-g(x,u)u}{h(|u|)} \geq
-h(x)
\end{eqnarray*}
with $\mu,h \in L^{\infty}(\Omega)$ and $\int_{\Omega} \mu (x) dx
> \int_{\partial \Omega} h(x) d\sigma $.

Let us state the energy functional.  Let $\Phi (u) =
-\int_{\Omega} F(x,u(x))dx, \quad \Gamma (u) = \int_{\partial
\Omega} G(x,u)d\sigma$ and $\psi (u) = \frac{1}{p} \nn Du
\nn_p^p$. Then our energy functional is $I = \psi+\Gamma + \Phi$
and is easy to check that is a $C^1$ functional and its critical
points are in fact weak solutions to problem (9).

\begin{lemma}
If Hypotheses $H(f,g),H(fg)$ holds, then $I:W^{1,p} \to \rr$
satisfies the $(PS)_c$ condition.
\end{lemma}

{\bf Proof:}

Let $\{ u_n \} \subseteq X$ be such that $|I(u_n)| \leq M  \in
\rr$ and
\begin{eqnarray*}
|< I^{'}(u_n), \phi>| \leq \vep_n \frac{\nn \phi \nn_{1,p}}{1+\nn
u_n \nn_{1,p}}, \mbox{ for all } \phi \in X, \quad \vep_n \to 0.
\end{eqnarray*}
Suppose that $u_n$ is unbounded. Then, at least for a subsequence,
we can say that $\nn u_n \nn_{1,p} \to \infty$. Let $y_n(x) =
\frac{u_n(x)}{\nn u_n \nn_{1,p}}$. Then it is easy to see that
$y_n \to y$ weakly in $X$ and strongly in $L^p(\Omega)$. From the
choice of the sequence we obtain
\begin{eqnarray}
| \frac{1}{p} \nn Du_n \nn_p^p  - \int_{\Omega} F(x,u_n(x))dx +
\int_{\partial \Omega} G(x,u_n(x))d\sigma| \leq M, \mbox{ for some
} M>0.
\end{eqnarray}
Dividing this inequality with $\nn u_n \nn^p_{1,p}$ we arrive at
\begin{eqnarray*}
| \frac{1}{p} \nn Dy_n \nn_p^p - \int_{\Omega}
\frac{F(x,u_n(x))}{\nn u_n \nn_{1,p}^p}dx + \int_{\partial \Omega}
\frac{G(x,u_n(x))}{\nn u_n \nn_{1,p}^p}d\sigma| \leq M.
\end{eqnarray*}

Using now $H(fg)$ we obtain $\nn Dy \nn_p^p = 0$. From this we
obtain that $y = \xi \neq 0$. Thus, we obtain $|u_n(x)| \to
\infty$.

Also we know that
\begin{eqnarray*}
|<I^{'}(u_n),\phi>| \leq \vep_n \frac{\nn \phi \nn_{1,p}}{1+\nn
u_n \nn_{1,p}}
\end{eqnarray*}
Choosing $\phi = u_n$ we arrive at
\begin{eqnarray}
|\nn Du_n \nn_p^p  - \int_{\Omega} f(x,u_n(x))u_n(x)dx +
\int_{\partial \Omega} g(x,u_n(x))u_n(x)d \sigma| \leq \vep_n
\frac{\nn u_n \nn_{1,p}}{1+\nn u_n \nn_{1,p}} .
\end{eqnarray}
Substituting (10) and (11) we obtain
\begin{eqnarray*}
& & -M - \vep_n \frac{\nn u_n \nn_{1,p}}{1+\nn u_n \nn_{1,p}}\leq
\\ & & \int_{\Omega} (p F(x,u_n(x)) - f(x,u_n(x))u_n(x))dx -
\int_{\partial \Omega} p G(x,u_n(x)) - g(x,u_n(x))u_n(x) d \sigma
\leq \\ & & M + \vep_n \frac{\nn u_n \nn_{1,p}}{1+\nn u_n
\nn_{1,p}}.
\end{eqnarray*}

Divide the last inequality with $h(\nn u_n \nn_{1,p})$ and using
$H(fg)$ we arrive at a contradiction as before. So $u_n$ is
bounded. In order to show that has a strongly convergent
subsequence we proceed by using well-known arguments (see
\cite{Halidias}).

\QED

\begin{lemma} If hypotheses $H(f,g), H(fg)$ holds, then there exists some
$e \in \rr$ such that $I(e) \leq 0$.
\end{lemma}

{\bf Proof:} In fact we are going to prove that there exists some
$a \in \rr$ big enough such that $I(a)\leq 0$. Suppose that this
is not the case. Then there exists a sequence $a_n \in \rr$ with
$a_n \to + \infty$ and $I(a_n ) \geq c > 0$. Using the same
arguments as before we can arrive at a contradiction.

\QED

\begin{lemma}If $H(f,g),H(fg)$ holds, then there exists some $\rho > 0$
such that for all $u \in W$ with $\nn u \nn_{1,p} = \rho$ we have
that $I(u) > \eta >0$.
\end{lemma}

{\bf Proof:} Suppose that this is not true. Then there exists a
sequence $\{ u_n \} \subseteq W$ such as $\nn u_n \nn_{1,p} =
\rho_n$ with $\rho_n \to 0$, with the property that $I(u_n) \leq
0$. So we arrive at
\begin{eqnarray*}
\nn Du_n \nn_p^p \leq p\int_{\Omega} F(x,u_n(x))dx- p
\int_{\partial \Omega} G(x,u_n(x))d \sigma.
\end{eqnarray*}
 Let $y_n(x) = \frac{u_n(x)}{\nn u_n \nn_{1,p}}$. We
can prove that there exists $\gamma > 0$ such that
\begin{eqnarray*}
p F(x,u) \leq (\theta (x) + \vep)|u|^p + \gamma |u|^{p^*}
\end{eqnarray*}
Take in account the last estimation and dividing with $\nn u_n
\nn_{1,p}^p$ we arrive at
\begin{eqnarray*}
\la_1 \nn y_n \nn_p^p \leq \nn Dy_n \nn_p^p  \leq \int_{\Omega}
(\theta (x) + \vep)|y_n(x)|^pdx + \gamma_1 \nn u_n \nn_p^{p^*-p}-
p \int_{\partial \Omega} \frac{G(x,u_n(x))}{\nn u_n \nn_{1,p}^p}
d\sigma.
\end{eqnarray*}
Recall that $y_n \to y$ strongly in $L^p(\Omega)$. Using the lower
semicontinuity of the norm we arrive at $\nn Dy \nn_p = \la_1 \nn
y \nn_p$. Note that $y_n \to y$ weakly in $X$ and recall that $\nn
Dy_n \nn_p \to \la_1 \nn y \nn_p = \nn Dy \nn_p$. Then from the
uniform convexity of $X$ we have $y_n \to y$ strongly  in $X$ and
$y \neq 0$. Thus $y \in W$ is an eigenfunction of $(-\Delta_p,
W)$. Then arguing as before we get the contradiction.

\QED

Then the existence of a nontrivial solution for problem (9)
follows from a variant of  Mountain-Pass Theorem (see Struwe
\cite{Struwe}, Thm. 8.4 and Example 8.2, or \cite{BBF}, Thm. 2.3
and Prop. 2.1). It is clear that we also extend the recently
results of Tang \cite{Tang} for Neumann problems.

\begin{remark}Take for example $f(u) = \frac{1}{u+1} + c u$, $c \in
\rr$ and $p=2$. Then it is easy to check that $2 F(u) - f(u)u =
ln(u+1)^2 - \frac{u}{u+1}$.

 It is easy to see that $h(u) = ln u$ satisfies all
the conditions and thus it is clear that the above results extend
all the known results (see for example \cite{BD}, \cite{Tang}).
\end{remark}

\end{document}